\documentclass%
[10pt,notitlepage]{article}
\usepackage{leftidx}
\usepackage{amsmath}
\usepackage{amssymb}
\usepackage{amsthm}
\usepackage{amscd}
\usepackage{latexsym}
\usepackage{graphicx}

\newtheorem{theorem}{Theorem}[section]
\newtheorem{proposition}[theorem]{Proposition}
\newtheorem{lemma}[theorem]{Lemma}

\newtheorem{definition}[theorem]{Definition}
\newtheorem{remark}[theorem]{Remark}

\def\Im{\operatorname{Im}}
\def\Ker{\operatorname{Ker}}

\def\Inn{\operatorname{Inn}}

\def\sign{\operatorname{sign}}
\def\Surj{\operatorname{Surj}}

\def\Diff{\operatorname{Diff}}
\def\Deck{\operatorname{Deck}}
\def\exp{\operatorname{exp}}
\def\mod{\operatorname{mod}}

\def\Sp{\operatorname{Sp}}

\def\Hom{\operatorname{Hom}}
\def\Prym{\operatorname{Prym}}
\begin{document}
\title{The abelianization of a symmetric mapping class group}
\author{Masatoshi Sato}
\date{}
\maketitle
\begin{abstract}
Let $\Sigma_{g,r}$ be a compact oriented surface of genus $g$ with $r$ boundary components. We determine the abelianization of the symmetric mapping class group $\hat{\mathcal{M}}_{(g,r)}(p_2)$ of a double unbranched cover $p_2:\Sigma_{2g-1,2r}\to\Sigma_{g,r}$ using the Riemann constant, Schottky theta constant, and the theta multiplier. We also give lower bounds of the abelianizations of some finite index subgroups of the mapping class group. 
\end{abstract}
\tableofcontents
\newpage
\section{Introduction}
Let $g$ be a positive integer, $r\ge0$, and $S$ a set of $n$ points in the interior of $\Sigma_{g,r}$. We denote by $\Diff_+(\Sigma_{g,r},\partial\Sigma_{g,r}, S)$ the group of all orientation preserving diffeomorphisms which fix the boundary $\partial\Sigma_{g,r}$ pointwise, and map $S$ onto itself allowing to permute it. The mapping class group $\mathcal{M}_{g,r}^n$ is the group of all isotopy classes $\pi_0\Diff_+(\Sigma_{g,r}, \partial\Sigma_{g,r}, S)$ of such diffeomorphisms. We write simply $\mathcal{M}_{g,r}:=\mathcal{M}_{g,r}^0$ and $\mathcal{M}_g^n:=\mathcal{M}_{g,0}^n$. The mapping class group and its finite index subgroups play an important role in low-dimensional topology, in the theory of Teichm\"{u}ller spaces and in algebraic geometry. For example, the level $d$ mapping class group $\mathcal{M}_{g,r}[d]$ is defined to be the finite index subgroup of $\mathcal{M}_{g,r}$ which acts trivially on $H_1(\Sigma_{g,r};\mathbf{Z}/d\mathbf{Z})$ for $d>0$. It arises as the orbifold fundamental group of the moduli space of genus $g$ curves with level $d$ structure.

To compute the abelianizations, or equivalently, the first integral homology groups of finite index subgroups is one of the important problems in the mapping class groups. The Torelli group $\mathcal{I}_{g,r}$ is the subgroup which acts trivially on $H_1(\Sigma_{g,r};\mathbf{Z})$. McCarthy\cite{mccarthy2000fcg} proved that the first rational homology group of a finite index subgroup that includes the Torelli group vanishes for $r=n=0$, and more generally, Hain\cite{hain28tga} proved it for any $r\ge0$, $n\ge0$.
\begin{theorem}[McCarthy, Hain]\label{theorem:finiteindex}
Let $\mathcal{M}$ be a finite index subgroup of $\mathcal{M}_{g,r}^n$ that includes the Torelli group where $g\ge3$, $r\geq0$. Then  
\[
H_1(\mathcal{M};\mathbf{Q})=0.
\]
\end{theorem}
This theorem gives us little information about $H_1(\mathcal{M};\mathbf{Z})$ as a finite group. In fact, Farb raised the problem to compute the abelianizations of the subgroup $\mathcal{M}_{g,r}[d]$ in \cite{farb2006spm} Problem 5.23 p.43.
In this paper, we confine ourselves to the case $r=0$ or $1$ when it is not specified. For a finite regular cover $p$ on $\Sigma_{g,r}$, possibly branched, Birman-Hilden\cite{birman1973ihr} defined the symmetric mapping class group $\hat{\mathcal{M}}_{(g,r)}(p)$. That is closely related to a finite index subgroup of the mapping class group. As stated in subsection \ref{subsection:defSMCG}, the symmetric mapping class group is a finite group extension of a certain finite index subgroup of the mapping class group. In particular, we will have $H_1(\hat{\mathcal{M}}_{(g,r)}(p);\mathbf{Q})=0$ for all abel covers $p$. But in general, the first integral homology groups of symmetric mapping class groups and finite index subgroups of $\mathcal{M}_{g,r}^n$ are unknown. 

One of the finite index subgroups, the spin mapping class group is defined by the subgroup of the mapping class group that preserves a spin structures on the surface. Lee-Miller-Weintraub\cite{lee1988rit} made the surjective homomophism from the spin mapping class group to $\mathbf{Z}/4\mathbf{Z}$ using the theta multiplier. Harer\cite{harer1993rpg} proved that this homomorphism is in fact an isomorphism.

In this paper, we determine the symmetric mapping class group $\hat{\mathcal{M}}_{(g,r)}(p_2)$ of an unbranched double cover $p_2:\Sigma_{2g-1,2r}\to\Sigma_{g,r}$ using the Riemann theta constant, Schottky theta constant, and the theta multiplier. We also compute a certain finite index subgroup $\mathcal{M}_{g,r}(p_2)$ of the mapping class group. That is included in the level 2 mapping class group $\mathcal{M}_{g,r}[2]$.

If we fix the symplectic basis of $H_1(\Sigma_{g,r};\mathbf{Z})$, the action of mapping class group $\mathcal{M}_{g,r}^n$ on $H_1(\Sigma_{g,r};\mathbf{Z})$ induces the surjective homomorphism
\[
\iota:\mathcal{M}_{g,r}\to Sp(2g;\mathbf{Z}),
\]
where $Sp(2g;\mathbf{Z})$ is the symplectic group of rank $2g$. Denote the image of $\mathcal{M}_{g,r}(p_2)$ under $\iota$ by $\Gamma_g(p_2)$. We also denote the image $\iota(\mathcal{M}_{g,r}[d])$ by $\Gamma_g[d]$, that is equal to the kernel $\Ker(Sp(2g;\mathbf{Z})\to Sp(2g;\mathbf{Z}/d\mathbf{Z}))$ of mod $d$ reduction. The main theorem is as follows. 
\begin{theorem}\label{main-theorem}
For $r=0,1$, when genus $g\ge4$,    
\begin{gather*}
H_1(\hat{\mathcal{M}}_{(g,r)}(p_2);\mathbf{Z})\cong H_1(\mathcal{M}_{g,1}(p_2);\mathbf{Z})\cong\mathbf{Z}/4\mathbf{Z},\\
H_1(\mathcal{M}_g(p_2);\mathbf{Z})\cong
\begin{cases}
\mathbf{Z}/4\mathbf{Z},\hspace{1cm} \text{ if\ \ } g:\text{odd},\\
\mathbf{Z}/2\mathbf{Z},\hspace{1cm} \text{ if\ \ } g:\text{even},
\end{cases}\\
H_1(\Gamma_g(p_2);\mathbf{Z})\cong\mathbf{Z}/2\mathbf{Z}.
\end{gather*}
\end{theorem}
After proving the theorem, we state that the first homology groups of the level $d$ mapping class group $H_1(\mathcal{M}_{g,1}[d];\mathbf{Z})$ have many elements of order 4 for any even integer $d$ (Proposition \ref{prop:leveld}).

In section \ref{symMCG}, we define the symmetric mapping class group, and describe the relation to a finite index subgroup of the mapping class group. In section \ref{genGamma}, we prove that the integral homology groups of $\hat{\mathcal{M}}_{(g,r)}(p_2)$ and $\mathcal{M}_{g,r}(p_2)$ are cyclic groups of order at most 4. We also have $H_1(\Gamma_g(p_2);\mathbf{Z})\cong \mathbf{Z}/2\mathbf{Z}$.
In section \ref{surj}, we construct an isomorphism $H_1(\hat{\mathcal{M}}_{(g,r)}(p_2);\mathbf{Z})\cong\mathbf{Z}/4\mathbf{Z}$ using the Schottky theta constant and the theta multiplier to complete the proof of theorem \ref{main-theorem}. 
\newpage
\section{The symmetric mapping class group}\label{symMCG}
In this section, we define the symmetric mapping class group following Birman-Hilden\cite{birman1973ihr}, and prove some properties. In particular, we describe $\mathcal{M}_{g,r}(p)=\Im P$ by means of the action of the mapping class group on the equivalent classes of the covers in Subsection \ref{subsection:cover}. We will see that the groups $\hat{\mathcal{M}}_{(g,r)}(p_2)$ and $\mathcal{M}_{g,r}(p_2)$ do not depend on the choice of the double cover $p_2$ up to isomorphism. 
\subsection{Definition of the symmetric mapping class group}\label{subsection:defSMCG}
Birman-Hilden\cite{birman1973ihr} defined the symmetric mapping class group of a regular cover $p: \Sigma_{g',r'}\to \Sigma_{g,r}$, possibly branched as follows. Denote the deck transformation group of the cover by $\Deck(p)$. 
\begin{definition}
Let $C(p)$ be the centralizer of the deck transformation $\Deck(p)$ in the diffeomorphism group $\Diff_+(\Sigma_{g',r'})$. The symmetric mapping class group of the cover $p$ is defined by 
\[
\hat{\mathcal{M}}_{(g,r)}(p)=\pi_0(C(p)\cap\Diff_+(\Sigma_{g',r'},\partial\Sigma_{g',r'})).
\] 
\end{definition}
Let $S\subset\Sigma_{g,r}$ be the branch set of the cover $p$. For $\hat{f}\in C(p)\cap \Diff_+(\Sigma_{g',r'},\partial\Sigma_{g',r'})$, there exists a unique diffeomorphism $f\in \Diff_+(\Sigma_{g,r},\partial\Sigma_{g,r}, S)$ such that the diagram
\[
\begin{CD}
\Sigma_{g',r'}@>\hat{f}>>\Sigma_{g',r'}\\
@V p VV @V p VV\\
\Sigma_{g,r}@>f>>\Sigma_{g,r}
\end{CD}
\]
commutes. Note that $f$ maps the branch set $S$ into itself. The diffeomorphism $f\in \Diff_+(\Sigma_{g,r},\partial\Sigma_{g,r}, S)$ is called the projection of $\hat{f}\in C(p)\cap \Diff_+(\Sigma_{g',r'},\partial\Sigma_{g',r'})$.   For $[\hat{f}], [\hat{g}]\in\hat{\mathcal{M}}_{(g,r)}(p)$ such that $[\hat{f}]=[\hat{g}]$, an isotopy between $\hat{f}$ and $\hat{g}$ induces the isotopy on the base space $\Sigma_{g,r}$ between the projections $f$ and $g$.  Hence we can define the homomorphism
\[
\begin{array}{cccc}
P:&\hat{\mathcal{M}}_{(g,r)}(p)&\to&\mathcal{M}_{g,r}^n,\\
&[\hat{f}]&\mapsto &[f]
\end{array}
\]
where $n\ge0$ is the order of $S$.
We denote the image $\Im P\subset\mathcal{M}_{g,r}^n$ by $\mathcal{M}_{g,r}(p)$. The kernel of $P$ is included in the group of isotopy classes of all the deck transformations in $\hat{\mathcal{M}}_{(g,r)}(p_2)$. Since any deck transformations without identity do not fix the boundary pointwise, we have $\Ker P=id$ when $r=1$. When $r=0$, $\Ker P$ consists of the isotopy classes of all the deck transformations. 
In particular, $\Ker P$ is a finite group. Apply the Lyndon-Hochschild-Serre spectral sequence to the group extension
\[1\to \Ker P\to \hat{\mathcal{M}}_{(g,r)}(p)\to\mathcal{M}_{g,r}(p)\to 0,\]
then we have 
\[
H_*(\hat{\mathcal{M}}_{(g,r)}(p);\mathbf{Q})\cong H_*(\mathcal{M}_{g,r}(p);\mathbf{Q}).
\]
\subsection{The action of the mapping class group on the equivalent classes of $G$-covers}\label{subsection:cover}
For a finite group $G$ and a finite set $S$, denote all the surjective homomorphisms $\pi_1(\Sigma_{g,r}-S, *)\to G$ by $\Surj(\pi_1(\Sigma_{g,r}-S, *), G)$. The group $G$ acts on this set by inner automorphism. Denote the quotient set by 
\[
m(G,*):=\Surj(\pi_1(\Sigma_{g,r}-S, *), G)/\Inn G.
\]
For paths $l,l':[0,1]\to\Sigma_{g,r}-S$ such that $l(0)=l'(1)$, we define $l\cdot l'$ to be the path obtained by traversing first $l'$ and then $l$. For a path $l:[0,1]\to\Sigma_{g,r}-S$, we define a isomorphism $l_*$ by
\[
\begin{array}{cccc}
l_*:&\pi_1(\Sigma_{g,r}-S,l(0))&\to&\pi_1(\Sigma_{g,r}-S,l(1)).\\
&\gamma&\mapsto&l\cdot\gamma\cdot l^{-1}
\end{array}
\]
If we pick a path $l$ from $*$ to $*'$, we have the isomorphism $l_*:\pi_1(\Sigma_{g,r}-S, *)\cong \pi_1(\Sigma_{g,r}-S, *')$.
Hence we also have the isomorphism
\[m(G,*)=m(G,*').\]
It is easy to see that this isomorphism does not depend on the choice of $l$, hence we denote $m(G):=m(G,*)$. 

The mapping class group $\mathcal{M}_{g,r}$ acts on the set $m(G)$. In fact the diffeomorphism $f\in \Diff_+(\Sigma_{g,r}, \partial\Sigma_{g,r}, S)$ induces the map
\[
\begin{array}{ccc}
m(G)&\to&m(G)\\
\empty[c]&\mapsto&[cf_*].
\end{array}
\]
\begin{proposition}
Let $c:\pi_1(\Sigma_{g,r}-S,*)\to G$ denote the monodromy homomorphism of a branched or unbranched $G$-cover $p:\Sigma_{g',r'}\to\Sigma_{g,r}$, where $S$ is the branch set. The stabilizer of $[c]\in m(G)$ is equal to $\mathcal{M}_{g,r}(p)$.
\end{proposition}
\begin{proof}
Suppose $[f]\in\mathcal{M}_{g,r}^n$ be in the stablizer of $[c]$. Since $[cf_*]=[c]$, there exists a path $l$ from $*$ to $f(*)$ such that
\[
c(l_*^{-1}f(\gamma))=c(\gamma), \text{ for } \gamma\in\pi_1(\Sigma_{g,r}-S,*).
\]
In particular, we have
\[
\Ker(c)=l_*^{-1}f_*(\Ker c).
\]
Hence the covers $p$ and $fp$ are equivalent. Choose a lift $\hat{l}$ of $l$, then there exists $\hat{f}\in\Diff(\Sigma_{g',r'})$ such that 
\[
p\hat{f}=fp:\Sigma_{g',r'}\to\Sigma_{g,r},\text{ and } \hat{f}(\hat{l}(0))=\hat{l}(1).
\]
Then we have
\[
\hat{f}c(\gamma)\hat{f}^{-1}=c(l_*^{-1}f(\gamma))=c(\gamma)\in\Diff_+\Sigma_{g',r'}.
\]
Hence $\hat{f}$ is in the centralizer $C(p)$ of the deck transformation group $\Deck(p)$. When $r=1$, $\Deck(p)$ acts on $\pi_0(\partial\Sigma_{g',r'})$ transitively. It is easy to see that for any $\hat{f}\in C(p)$, there exist $t\in \Deck(p)$ such that $\hat{f}t$ acts trivially 
on $\pi_0(\partial\Sigma_{g',r'})$. Therefore, there exists $t\in \Deck(p)$ such that $\hat{f}t\in C(p)\cap\Diff_+(\Sigma_{g',r'},\partial\Sigma_{g',r'})$ and $f=P([\hat{f}t])$.

Conversely, suppose $f=P(\hat{f})\in\mathcal{M}_{g,r}(p)$. Choose a path $\hat{l}$ such that $\hat{f}(\hat{l}(0))=\hat{l}(1)$. Denote the projection $l=p\hat{l}$, then we have
\[
c(l_*^{-1}f(\gamma))=\hat{f}c(\gamma)\hat{f}^{-1}=c(\gamma)\in\Diff_+\Sigma_{g',r'}.
\]
Hence we have $[c]=[cf_*]$.
\end{proof}
Hence, $\mathcal{M}_{g,r}(p)$ is a finite index subgroup of the mapping class group. In particular, if $p$ is an abel cover, $\mathcal{M}_{g,r}(p)$ includes the Torelli group. By Theorem \ref{theorem:finiteindex}, we have $H_1(\mathcal{M}_{g,r}(p);\mathbf{Q})=0$. Consider the double covers on $\Sigma_{g,r}$. The number of the equivalent classes of double unbranched covers on $\Sigma_{g,r}$ are $2^{2g}-1$. Since the action of mapping class group $\mathcal{M}_{g,r}$ on $m(\mathbf{Z}/2\mathbf{Z})$ is transitive, the subgroup $\mathcal{M}_{g,r}(p_2)$ does not depend on the choice of the double cover $p_2$ up to conjugate. It is easy to see that $\hat{\mathcal{M}}_{(g,r)}(p_2)$ is also unique up to isomorphism.
\newpage
\section{A lower bound of the order of the cyclic group $H_1(\hat{\mathcal{M}}_{(g,r)}(p_2);\mathbf{Z})$}\label{genGamma}
In this section we prove that the integral homology groups of $\hat{\mathcal{M}}_{(g,r)}(p_2)$ and $\mathcal{M}_{g,r}(p_2)$ are cyclic groups of order at most $4$. We compute $H_1(\Gamma_g(p_2);\mathbf{Z})$ in Subsection \ref{subsection:Gamma} and $H_1(\mathcal{I}_{g,r};\mathbf{Z})_{\mathcal{M}_{g,r}(p_2)}$ in Subsection \ref{subsection:torelli} to obtain the lower bound.

In subsection \ref{subsection:cover}, we proved that the symmetric mapping class group $\hat{\mathcal{M}}_{(g,r)}(p_2)$ and $\mathcal{M}_{g,r}(p_2)$ do not depend on the choice of the unbranched double cover $p_2$ up to isomorphism. Hence we fix the unbranched double cover $p_2$ whose monodromy $c\in\Hom(\pi_1(\Sigma_{g,r});\mathbf{Z}/2\mathbf{Z})\cong H^1(\Sigma_{g,r};\mathbf{Z}/2\mathbf{Z})$ is equal to the Poincar\'e dual of $B_g$ in Figure \ref{symplectic}.
\subsection{The first homology group $H_1(\Gamma_g(p_2);\mathbf{Z})$}\label{subsection:Gamma}
In this subsection, using the generators of $\Gamma_g[2]$ in Igusa\cite{igusa1964grt}, we prove that $H_1(\Gamma_g(p_2);\mathbf{Z})$ is a cyclic group of order 2. We also prove that $H_1(\hat{\mathcal{M}}_{(g,r)}(p_2);\mathbf{Z})$ and $H_1(\mathcal{M}_{g,r}(p_2);\mathbf{Z})$ are cyclic of order at most 4 when genus $g\ge4$, using the $\mathcal{M}_{g,r}$ module structure  of the abelianization of the Torelli group determined by Johnson\cite{johnson1980aqm}. In particular, we obtain $H_1(\mathcal{M}_g(p_2);\mathbf{Z})\cong \mathbf{Z}/2\mathbf{Z}$ if genus $g\ge 4$ is even. In the next section, we complete the proof of Theorem \ref{main-theorem}. 

\begin{figure}[htbp]
\begin{center}
\includegraphics{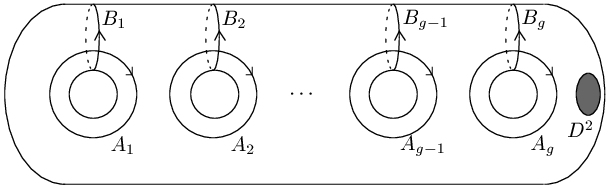}
\end{center}
\caption{}
\label{symplectic}
\end{figure}

\begin{figure}[htbp]
  \begin{center}
    \includegraphics{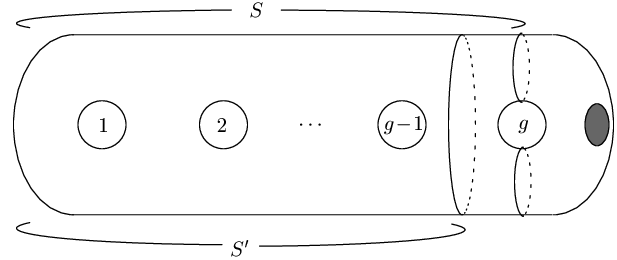}
  \end{center}
  \caption{}
  \label{subsurface}
\end{figure}

We consider $\Sigma_{g,1}=\Sigma_g-D^2\subset\Sigma_g$. 
Pick simple closed curves $\{A_i, B_i\}_{i=1}^g\subset \Sigma_{g,r}$ as shown in Figure \ref{symplectic}. They give a symplectic basis of $H:=H_1(\Sigma_{g,r};\mathbf{Z})$ which we denote by the same symbol $\{A_i, B_i\}_{i=1}^g$. The action of the mapping class group on $H_1(\Sigma_{g,r};\mathbf{Z})$ induces 
\[
\iota:\mathcal{M}_{g,r}\to\Sp(2g,\mathbf{Z}).
\]
We denote the Dehn twist along the simple closed curve $A_g$ by $a\in\mathcal{M}_{g,r}$. 
Let $S$ be a subsurface in $\Sigma_{g,r}$ as shown in Figure \ref{subsurface} and denote their mapping class groups which fix the boundary pointwise by $\mathcal{M}_S$. 
The inclusion $S\to \Sigma_g$ induces a homomorphism
\[
i_S:\mathcal{M}_S\to \mathcal{M}_g.
\]
As in Introduction, we denote by $\iota: \mathcal{M}_g\to\Sp(2g;\mathbf{Z})$ the homomorphism defined by the action of $\mathcal{M}_g$ on the homology group $H$, and denote the ring of integral $n$-square matrices by $M(n;\mathbf{Z})$ for a positive integer $n$. It is easy to see that the image of $i_S(\mathcal{M}_S)$ under $\iota$ is
\[
\iota(i_S(\mathcal{M}_S))=
\left\{
\sigma=\left.
\begin{pmatrix}
\alpha'&\leftidx{^t}{v_1}{}&\beta'&0\\
0&1&0&0\\
\gamma'&\leftidx{^t}{v_2}{}&\delta'&0\\
v_3&k&v_4&1
\end{pmatrix}
\in \Sp(2g;\mathbf{Z})
\ \right|\ 
\begin{array}{c}
\alpha', \beta', \gamma', \delta'\in M(g-1;\mathbf{Z}),\\
v_1,v_2,v_3,v_4\in \mathbf{Z}^{g-1}, k\in\mathbf{Z}
\end{array}
\right\}.
\]
\begin{proposition}\label{generator}
When $g\ge1$, $\Gamma_g(p_2)$ is generated by $\iota(i_S(\mathcal{M}_S))$ and $ \iota(a^2)$.
\end{proposition}
\begin{proof}
First, we show that $\Gamma_g(p_2)$ is generated by $\iota(i_S(\mathcal{M}_S))$ and $\Gamma_g[2]$.  Since an element $\sigma\in\Gamma_g(p_2)$ preserves the homology class $B_g\in H_1(\Sigma_{g,r};\mathbf{Z}/2\mathbf{Z})$, it can be written in the form
\[
\sigma\equiv
\begin{pmatrix}
\alpha'&\leftidx{^t}{v_1}{}&\beta'&0\\
0&1&0&0\\
\gamma'&\leftidx{^t}{v_2}{}&\delta'&0\\
v_3&k&v_4&1
\end{pmatrix}
\mod 2.\]
Hence there exists $\sigma_0\in\iota(i_S(\mathcal{M}_S))$ such that 
\[\sigma_0\equiv \sigma \ \mod 2, \]
so that $\Gamma_g(p_2)$ is generated by $\iota(i_S(\mathcal{M}_S))$ and $\Gamma_g[2]$.

Next, we describe the generators of $\Gamma_g[2]$ given in Igusa\cite{igusa1964grt}. We denote by $I_n$ the unit matrix of order $n$, and by $e_{ij}$ the $2g$-square matrix with 1 at the $(i,j)$-th entry and 0 elsewhere. As was shown in Igusa\cite{igusa1964grt}, $\Gamma_g[2]$ is generated by
\begin{align*}
\alpha_{ij}=&I_{2g}+2e_{ij}-2e_{g+j,g+i}&& 1\le i,j\le g,\  i\ne j, \\
\alpha_{ii}=&I_{2g}-2e_{ii}-2e_{i+g,i+g}&& 1\le i\le g, \\
\beta_{ij}=&I_{2g}+2e_{i,j+g}+2e_{j,i+g}&& 1\le i<j\le g, \\
\beta_{ii}=&I_{2g}+2e_{i,i+g}&& 1\le i\le g, \\
\gamma_{ij}=&\leftidx{^t}{\beta_{ij}}{}&&1\le i\le j\le g.
\end{align*}
To prove the proposition, it suffices to show that these matrices are in the subgroup of $\Gamma_g(p_2)$ generated by $\iota(i_S(\mathcal{M}_S))$ and $\iota(a^2)$. The matrices
\[
\alpha_{ij} (1\le i\le g-1, 1\le j\le g),\  
\beta_{ij} (1\le i\le j\le g-1)\text{, and}\  
\gamma_{ij} (1\le i\le j\le g)
\]
are clearly in $\iota(i_S(\mathcal{M}_S))$. Choose oriented simple closed curves $C_i, C'_i, C_{ij}, C'_{ij}, C''_{ij}\subset \Sigma_{g,r}$ such that $[C_i]=A_i$, $[C'_i]=B_i$, $[C_{ij}]=A_i+A_j$, $[C'_{ij}]=B_i+B_j$, $[C''_{ij}]=A_i+B_j$. Denote the Dehn twist along a simple closed curve $C$ by $T_C$. Then the matrices
\[
\alpha_{gj} (1\le j\le g-1) \ \text{ and }\ \beta_{ig} (1\le i\le g-1)
\]
are written as $\iota(T_{C''_{gi}}^2T_{C'_i}^{-2}T_{C_g}^{-2})$ and $\iota(T_{C_{ig}}^2T_{C_i}^{-2}T_{C_g}^{-2})$ respectively.  Clearly $\iota(T_{C'_i}^2)$ and $\iota(T_{C_i}^2)$ are in $\iota(i_S(\mathcal{M}_S))$, and we have $\iota(T_{C_g}^2)=\iota(a^2)$.  Denote the two boundary components of $S$ by $S_1$ and $S_2$. For any two arcs $l_1,l_2:[0,1]\to S$ that satisfy $l_1(0)=l_2(0)\in S_1$ and $l_1(1)=l_2(1)\in S_2$, there exists $\varphi\in \mathcal{M}_S$ such that
\[\varphi l_1=l_2.\]
Choose $C''_{gi}$ and $C_{ig}$ such that $\sharp(C''_{gi}\cap S_1)=\sharp(C_{ig}\cap S_1)=1$ and they intersect with $S_1$ transversely, there exist $\psi, \psi'\in i_S(\mathcal{M}_S)$ that satisfy $[\psi(C''_{gi})]=[\psi'({C_{ig}})]=A_g$. Thus we have  
\[\psi T_{C''_{gi}}^2\psi^{-1}=\psi'T_{C_{ig}}^2{\psi'}^{-1}=a^2.\]
This proves the matrices $\alpha_{gj}$ and $\beta_{ig}$ are in the subgroup of $\Gamma_g(p_2)$ generated by $\iota(i_S(\mathcal{M}_S))$ and $\iota(a^2)$. Finally the matrices $\alpha_{gg}$ and $\beta_{gg}$ satisfy $\alpha_{gg}=\iota(T_{C''_{gg}}^2)\beta_{gg}\gamma_{gg}^{-1}$, and $\beta_{gg}=\iota(a^2)$. Hence $\alpha_{gg}$ and $\beta_{gg}$ are also in the subgroup, as was to be shown.
\end{proof}
Using Proposition \ref{generator}, we now calculate the first homology group $H_1(\Gamma_g(p_2);\mathbf{Z})$. 
\begin{proposition}\label{Gamma(p_2)}
When $g\ge4$, 
\[
H_1(\Gamma_g(p_2);\mathbf{Z})\cong \mathbf{Z}/2\mathbf{Z}.
\] 
\end{proposition}
\begin{proof}
Powell \cite{powell1978ttm} had proved $H_1(\mathcal{M}_g;\mathbf{Z})=0$ when $g\ge3$. More generally, Harer \cite{harer1983shg} proved that $H_1(\mathcal{M}_{g,r};\mathbf{Z})=0$ when $g\ge3$ for any $r$. Hence the first homology $H_1(\mathcal{M}_S;\mathbf{Z})$ vanishes since genus of $S$ $\ge3$. We have
\[i_S(\mathcal{M}_S)=\{0\}\subset H_1(\mathcal{M}_g(p_2);\mathbf{Z}).\]
Since we proved that the group $\Gamma_g(p_2)$ is generated by $\iota(i_S(\mathcal{M}_S))$ and $\iota(a^2)$ in Proposition \ref{generator}, the homology group $H_1(\Gamma_g(p_2);\mathbf{Z})$ is generated by $[\iota(a^2)]$. 

Next, we construct a surjective homomorphism $\Gamma_g(p_2)\to\mathbf{Z}/2\mathbf{Z}$. Since any $\sigma=(\sigma_{ij})\in\Gamma_g(p_2)$ preserves the homology class $B_g\in H_1(\Sigma_{g,r};\mathbf{Z}/2\mathbf{Z})$, we have
\[\sigma_{gi}\equiv \delta_{ig}\text{, and \ } \sigma_{i\,2g}\equiv \delta_{i\,2g}\ \mod 2,\]
where $\delta$ is the Kronecker delta. Then for $\sigma, \sigma'\in\Gamma_g(p_2)$, the $(g,2g)$-th entry of $\sigma\sigma'$ satisfies
\[(\sigma\sigma')_{g\,2g}=\sum_{i=1}^{2g}\sigma_{gi}\sigma'_{i\,2g}\equiv \sigma_{g\,2g}+\sigma'_{g\,2g}\ \mod 4.\]
Hence we have the homomorphism
\[
\begin{array}{cccc}
\Psi:&\Gamma_g(p_2)&\to&\mathbf{Z}/2\mathbf{Z}\\
&\sigma&\mapsto& \displaystyle\frac{\sigma_{g\,2g}}{2}. 
\end{array}
\]
Since $\Psi([\iota(a^2)])=1$, we have $[\iota(a^2)]\ne0\in H_1(\Gamma_g(p_2);\mathbf{Z})$.  

Finally, to complete the proof it suffices to show that $2[\iota(a^2)]=0$. Apply the Lyndon-Hochschild-Serre spectral sequence to the group extension
\[1\to \mathcal{I}_{g,r}\to \mathcal{M}_{g,r}(p_2)\to\Gamma_g(p_2)\to 0,\]
then we have 
\[
H_1(\mathcal{I}_{g,r};\mathbf{Z})_{\mathcal{M}_{g,r}(p_2)}\to H_1(\mathcal{M}_{g,r}(p_2);\mathbf{Z})\to H_1(\Gamma_g(p_2);\mathbf{Z})\to 0.
\]
Denote by $D$ and $D'$ the simple closed curves as shown in Figure \ref{t_dt_d'}. 
\begin{figure}[htbp]
  \begin{center}
    \includegraphics{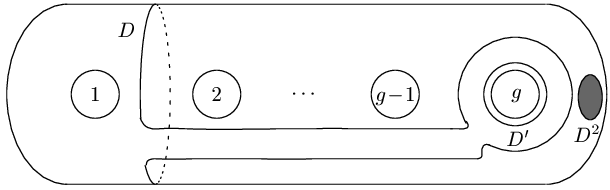}
  \end{center}
  \caption{}
  \label{t_dt_d'}
\end{figure}
Denote by $c_1$, $c_2$, and $c_3$ the Dehn twists along the simple closed curves $C_1$, $C_2$, and $C_3$ as shown in Figure \ref{abc} respectively. Since $c_1$ and $c_2$ are in $i_S(\mathcal{M}_S)$, $[c_1]=[c_2]=0\in H_1(\mathcal{M}_{g,r}(p_2);\mathbf{Z})$. By the chain relation, we have $T_DT_{D'}=(c_1c_2c_3)^4$. 
\begin{figure}[htbp]
  \begin{center}
    \includegraphics{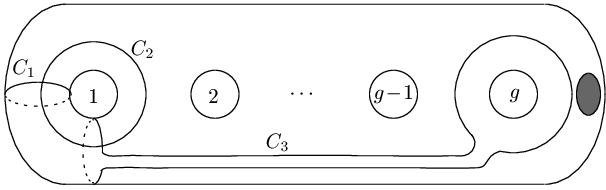}
  \end{center}
  \caption{}
  \label{abc}
\end{figure}
Using the braid relations $c_1c_3=c_3c_1$, $c_1c_2c_1=c_2c_1c_2$, and $c_2c_3c_2=c_3c_2c_3$, we have
\[
[T_DT_{D'}]=[(c_1c_2c_3)^4]=[c_3c_2c_1^2c_2c_3]=[c_3c_2c_1^2c_2^{-1}c_3^{-1}]+[c_3c_2^2c_3^{-1}]+[c_3^2]
\in H_1(\mathcal{M}_{g,r}(p_2);\mathbf{Z}).\]
Since $c_3c_2c_1^2c_2^{-1}c_3^{-1}$ and $c_2c_1^2c_2^{-1}$ are the squares of the Dehn twists along the simple closed curves $c_3c_2(C_1)$ and $c_3(C_2)$, we have
\[
[c_3c_2c_1^2c_2^{-1}c_3^{-1}]=[c_3c_2^2c_3^{-1}]=[c_3^2]=[a^2].
\]
Hence $[T_DT_{D'}^{-1}]=[T_DT_{D'}]+[T_{D'}^{-2}]=2[a^2]$. Since $T_DT_{D'}^{-1}\in\mathcal{I}_{g,r}$, it follows that $2[\iota(a^2)]=[\iota(T_DT_{D'}^{-1})]=0\in H_1(\Gamma_g(p_2);\mathbf{Z})$. 
This proves the proposition. 
\end{proof}
\subsection{The coinvariant $H_1(\mathcal{I}_{g,r};\mathbf{Z})_{\mathcal{M}_{g,r}(p_2)}$}\label{subsection:torelli}
To calculate the first homology group of the symmetric mapping class groups, we compute $H_1(\mathcal{I}_{g,r};\mathbf{Z})_{\mathcal{M}_{g,r}(p_2)}$. 
\begin{lemma}\label{torelli}
When $g\ge4$,  
\[
H_1(\mathcal{I}_g;\mathbf{Z})_{\mathcal{M}_g(p_2)}\cong
\begin{cases}
\mathbf{Z}/2\mathbf{Z},\hspace{0.3cm}&\text{if }g:\text{odd},\\
0,\hspace{0.3cm}&\text{if }g:\text{even},
\end{cases}
\]
\[
H_1(\mathcal{I}_{g,1};\mathbf{Z})_{\mathcal{M}_{g,1}(p_2)}\cong\mathbf{Z}/2\mathbf{Z}.
\]
Moreover $H_1(\mathcal{I}_{g,r};\mathbf{Z})_{\mathcal{M}_{g,r}(p_2)}$ is generated by $T_DT_{D'}^{-1}\in\mathcal{I}_{g,r}$ for $r=0,1$. 
\end{lemma}
Before proving the lemma, we review the space of boolean polynomials. Let $H$ denote the first homology group $H_1(\Sigma_{g,r};\mathbf{Z})$ of the surface as before. Consider the polynomial ring with coefficients in $\mathbf{Z}/2\mathbf{Z}$ with the basis $\bar{x}$ for $x\in H\otimes \mathbf{Z}/2\mathbf{Z}$. Denote by $J$ the ideal in the polynomial generated by 
\[\overline{x+y}-(\bar{x}+\bar{y}+x\cdot y),\hspace{0.4cm}\bar{x}^2-\bar{x},\hspace{0.4cm}\text{ for } x,y\in H\otimes \mathbf{Z}/2\mathbf{Z}.\]
The space of boolean polynomials of degree at most $n$ is defined by
\[B^n=\frac{M_n}{J\cap M_n},\]
where $M_n$ is the module of all polynomials of degree at most $n$. Note that $B^n$ is isomorphic to the $\mathbf{Z}/2\mathbf{Z}$ module of all square free polynimials of degree at most $n$ generated by $\{\bar{A}_i, \bar{B}_i\}_{i=1}^g$.

Denote $B^3$ by $B_{g,1}^3$, and for $\alpha=\Sigma_{i=1}^g\bar{A}_i\bar{B}_i\in B^2$,  the cokernel of
\[
\begin{array}{ccc}
B^1&\to&B^3\\
x&\mapsto&\alpha x
\end{array}
\]
by $B_{g,0}^3$. The action of $\mathcal{M}_{g,r}$ on $H$ induces an action on $B_{g,r}^3$. Birman-Craggs\cite{birman1978muim} defined a family of homomorphisms $\mathcal{I}_g\to \mathbf{Z}/2\mathbf{Z}$. Johnson\cite{johnson1980qfa} showed that these homomorphisms give a surjective homomorphism of $\mathcal{M}_{g,r}$ modules
\[
\mu: \mathcal{I}_{g,r}\to B_{g,r}^3.
\]
For $r=0,1$, Johnson\cite{johnson1985stg} showed that the induced homomorphism $H_1(\mathcal{I}_{g,r};\mathbf{Z})_{\mathcal{M}_{g,r}[2]}\cong B_{g,r}^3$ is an isomorphism. 
\begin{proof}[proof of Lemma \ref{torelli}]
Since $\mu$ is an isomorphism of $\mathcal{M}_{g,r}$ module, we have  
\[
H_1(\mathcal{I}_{g,r};\mathbf{Z})_{\mathcal{M}_{g,r}(p_2)}\cong (B_{g,r}^3)_{\mathcal{M}_{g,r}(p_2)}.
\]
Hence it suffices to compute $(B_{g,r}^3)_{\mathcal{M}_{g,r}(p_2)}$ to prove the lemma. Denote the subsurface $S'\subset S$ of genus $g-1$ as shown in Figure \ref{subsurface}. $\mathcal{I}_{S'}$ is the Torelli group of $S'$, that is the subgroup of $\mathcal{M}_{S'}$ which act trivially on $H_1(S';\mathbf{Z})$. 
Consider the homomorphism 
\[
(\mathcal{I}_{S'})_{\mathcal{M}_{S'}}\to (\mathcal{I}_{g,r})_{\mathcal{M}_{g,r}(p_2)}\cong (B_{g,r}^3)_{\mathcal{M}_{g,r}(p_2)}.
\]
induced by the inclusion $S'\to\Sigma_{g,r}$.
Since $(\mathcal{I}_{S'})_{\mathcal{M}_{S'}}=0$ (Johnson\cite{johnson1979hsa}), the image of the homomorphism is trivial. Thus we have  
\[\bar{1}=\bar{X}=\bar{X}\bar{Y}=\bar{X}\bar{Y}\bar{Z}=0, \ \text{ for } \{X,Y,Z\}\subset\{A_1,A_2,\cdots,A_{g-1},B_1,B_2\cdots,B_{g-1}\}.\] 

For $X=A_g, B_g$, we have
\begin{gather*}
(I_{2g}+e_{1,g+1})(\bar{B}_1\bar{X})=(\bar{B}_1+\bar{A}_1+1)\bar{X},\ 
(I_{2g}+e_{g+1,1})(\bar{A}_1\bar{X})=(\bar{A}_1+\bar{B}_1+1)\bar{X},\\
\text{ and }\  (I_{2g}+e_{1,2}-e_{g+2,g+1})(\bar{A}_2\bar{X})=(\bar{A}_2+\bar{A}_1)\bar{X}.
\end{gather*}
Hence $\bar{X}=\bar{A}_1\bar{X}=\bar{B}_1\bar{X}=0\in(B_{g,r}^3)_{\mathcal{M}_{g,r}(p_2)}$.  For $1< i< g$, we have
\[(I_{2g}+e_{g+i,1}+e_{g+1,i})(\bar{A}_1\bar{X})=(\bar{A}_1+\bar{B}_i)\bar{X}
,\text{ and }\ (I_{2g}+e_{i,g+1}+e_{1,g+i})(\bar{B}_1\bar{X})=(\bar{B}_1+\bar{A}_i)\bar{X}.\] 
Hence $\bar{B}_i\bar{X}=\bar{A}_i\bar{X}=0\in(B_{g,r}^3)_{\mathcal{M}_{g,r}(p_2)}$. If we put $\bar{X}=\bar{A}_g\bar{B}_g$, we have $\bar{Y}\bar{A}_g\bar{B}_g=0$ in the same way for $\bar{Y}\in\{1,\bar{A}_1,\bar{A}_2,\cdots,\bar{A}_{g-1},\bar{B}_1,\bar{B}_2,\cdots,\bar{B}_{g-1}\}$. For $X=A_g, B_g$, and any $i,j$ such that $1\le i,j<g$, $i\ne j$, we have
\begin{gather*}
(I_{2g}+e_{g+j,j})(\bar{A}_i\bar{A}_j\bar{X})=\bar{A}_i\bar{A}_j\bar{X}+\bar{A}_i\bar{B}_j\bar{X}+\bar{A}_i\bar{X}, \ (I_{2g}+e_{j,g+j})(\bar{A}_i\bar{B}_j\bar{X})=\bar{A}_i\bar{B}_j\bar{X}+\bar{A_i}\bar{A}_j\bar{X}+\bar{A}_i\bar{X},\\
(I_{2g}+e_{g+j,j})(\bar{B}_i\bar{A}_j\bar{X})=\bar{B}_i\bar{A}_j\bar{X}+\bar{B}_i\bar{B}_j\bar{X}+\bar{B}_i\bar{X},\ 
(I_{2g}+e_{g+i,g}+e_{2g,g+i})(\bar{A}_i\bar{A}_g\bar{B}_g)=\bar{A}_i\bar{A}_g\bar{B}_g+\bar{A}_i\Bar{B}_i\bar{B}_g,\\
\text{and }\ (I_{2g}-e_{1,1}-e_{g+1,g+1}+e_{i,1}+e_{1,i}+e_{g+i,g+1}+e_{g+1,g+i})(\bar{A}_1\bar{B}_1\bar{A}_g)=\bar{A}_i\bar{B}_i\bar{A}_g.
\end{gather*}
Hence $\bar{A}_i\bar{B}_j\bar{X}=\bar{A}_i\bar{A}_j\bar{X}=\bar{B}_i\bar{B}_j\bar{X}=\bar{A}_i\bar{B}_i\bar{B}_g=0$, and $\bar{A}_1\bar{B}_1\bar{A}_g=\bar{A}_i\bar{B}_i\bar{A}_g\in (B_{g,r}^3)_{\mathcal{M}_{g,r}(p_2)}$.

Therefore $(B_{g,r}^3)_{\mathcal{M}_{g,r}(p_2)}$ is a cyclic group of order 2 with generator $\bar{A}_1\bar{B}_1\bar{A}_g$ or a trivial group.  For $r=0$, $B_{g,0}^3$ has a relation
\[
\alpha\bar{A}_g=(\sum_{i=1}^g\bar{A}_i\bar{B}_i)\bar{A}_g=0,
\]
so that $(g-1)\bar{A}_1\bar{B}_1\bar{A}_g=0\in(B_{g,0}^3)_{\mathcal{M}_g(p_2)}$. This shows that $(B_{g,r}^3)_{\mathcal{M}_{g,r}(p_2)}$ is trivial when $g$ is even and $r=0$.

Next we consider the case $g$ is odd or $r=1$. Let $S_n$ be the permutation group of degree $n$ and $\sign(s)$ the sign of $s\in S_n$. Denote by $\Lambda^nH$ the image of the homomorphism 
\[
\begin{array}{cccc}
\lambda:&H^{\otimes n}&\to&H^{\otimes n}\\
&x_1\otimes x_2\otimes\cdots\otimes x_n&\mapsto&\displaystyle\sum_{s\in S_n}\sign(s) x_{s(1)}\otimes x_{s(2)}\otimes\cdots\otimes x_{s(n)}.
\end{array}
\]
Denote by $V_1$ and $V_0$ the module $\Lambda^3H$ and the cokernel of 
\[
\begin{array}{ccc}
H&\to&\Lambda^3H\\
X&\mapsto&\sum_{i=1}^g A_i\wedge B_i\wedge X,
\end{array}
\]
respectively. Then Johnson\cite{johnson1980aqm} shows 
\[
\begin{array}{cccc}
\displaystyle\frac{B_{g,r}^3}{B^2}&\to &V_r\otimes\mathbf{Z}/2\mathbf{Z}\\[3pt]
\bar{X}\bar{Y}\bar{Z}&\mapsto &X\wedge Y\wedge Z,
\end{array}
\]
is a well-defined $\mathcal{M}_{g,r}$ module isomorphism. Now we have a $\mathcal{M}_{g,r}(p_2)$ homomorphism 
\[
\begin{array}{cccc}
(B_g\cdot)C:&(B_{g,r}^3)_{\mathcal{M}_{g,r}(p_2)}&\to &\mathbf{Z}/2\mathbf{Z}\\[3pt]
&\bar{X}\bar{Y}\bar{Z}&\mapsto &(X\cdot Y)B_g\cdot Z+(Y\cdot Z)B_g\cdot X+(Z\cdot X)B_g\cdot Y.
\end{array}
\]
Here it should be remarked that the intersection number with $B_g$ $(B_g\cdot):H\otimes\mathbf{Z}/2\mathbf{Z}\to\mathbf{Z}/2\mathbf{Z}$ is $\mathcal{M}_{g,r}(p_2)$-invariant. Since $(B_g\cdot)C (\bar{A}_1\bar{B}_1\bar{A}_g)=1$, it is surjective. Hence $(B_{g,r}^3)_{\mathcal{M}_{g,r}(p_2)}$ is a cyclic group of order 2 with generator $\bar{A}_1\bar{B}_1\bar{A}_g$. Johnson\cite{johnson1980qfa} computed $\mu(T_DT_{D'}^{-1})=\bar{A}_1\bar{B}_1(\bar{A}_g+1)$, so that $T_DT_{D'}^{-1}$ is a generator of $H_1(\mathcal{I}_{g,r};\mathbf{Z})_{\mathcal{M}_{g,r}(p_2)}$.
\end{proof}
Now, we prove that $H_1(\hat{\mathcal{M}}_{(g,r)};\mathbf{Z})$ and $H_1(\mathcal{M}_{g,r};\mathbf{Z})$ are cyclic groups of order at most 4.
We need the following Lemma. 
\begin{lemma}
Let $\hat{b}:\hat{\mathcal{M}}_{(g,1)}(p_2)\to\hat{\mathcal{M}}_{(g)}(p_2)$ be a homomorphism induced by an obvious embedding $\Sigma_{2g-1,2}\to\Sigma_{2g-1}$. Then $\hat{b}$ is surjective.
\end{lemma}
\begin{proof}
By the obvious embedding $\Sigma_{g,1}\to\Sigma_{g}$, we have a surjective homomorphism $b:\mathcal{M}_{g,1}(p_2)\to\mathcal{M}_g(p_2)$. Since the diagram
\[
\begin{CD}
\hat{\mathcal{M}}_{(g,1)}(p_2)@>\hat{b}>>\hat{\mathcal{M}}_{(g)}(p_2)\\
@V P VV @V P VV\\
\mathcal{M}_{g,1}(p_2)@>b>>\mathcal{M}_{g}(p_2)
\end{CD}
\]
commutes, $\hat{b}P=Pb$ is surjective. Hence it surfices to show $\Ker P\subset\hat{\mathcal{M}}_{(g)}(p_2)$ is included in $\Im\hat{b}$. Recall that $\Ker P$ consists of the isotopy classes of all the deck transformation.
% }'Ì'}"ü
\begin{figure}[htbp]
  \begin{center}
    \includegraphics{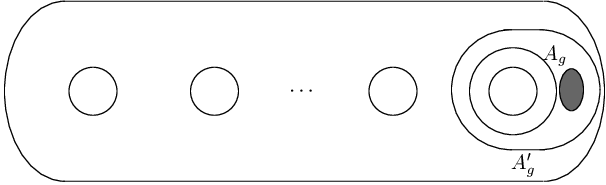}
  \end{center}
  \caption{}
  \label{fig:a_g}
\end{figure}
Cut the surface $\Sigma_{g,r}$ along the two simple closed curves $A_g$, $A'_g$ in Figure \ref{fig:a_g}. Then we have the subsurface $S_0$ of genus $g-1$ and the other subsurface $S'_0$ of genus $0$. We can construct a diffeomorphism $\hat{f}_0\in C(p)\cap\Diff_+(\Sigma_{2g-1,2},\partial\Sigma_{2g-1,2})$ which have the following properties:
\begin{enumerate}
\item $\hat{f}_0|_{p^{-1}(S_0)}$ is the restriction of the deck transformation $t\ne id$,
\item $\hat{f}_0|_{p^{-1}(S_1)}={T'}_{A_g}{T'}_{A_g'}^{-1}$, where ${T'}_{A_g}$ and ${T'}_{A_g'}$ is the half Dehn twists along $A_g$ and $A_g'$.
\end{enumerate}
Then $\hat{f}_0$ is included in $C(p)\cap\Diff_+(\Sigma_{2g-1,2},\partial\Sigma_{2g-1,2})$, and the image of $[\hat{f}_0]$ under $\hat{\mathcal{M}}_{(g,1)}(p_2)\to\hat{\mathcal{M}}_{(g)}(p_2)$ equals the deck transformation $t$. This proves the lemma.
\end{proof}
In the proof of Proposition \ref{Gamma(p_2)}, we have the exact sequence
\[
H_1(\mathcal{I}_{g,r};\mathbf{Z})_{\mathcal{M}_{g,r}(p_2)}\to H_1(\mathcal{M}_{g,r}(p_2);\mathbf{Z})\to H_1(\Gamma_g(p_2);\mathbf{Z})\to 0.
\]
By Proposition \ref{Gamma(p_2)} and Lemma \ref{torelli}, we obtain
\[
H_1(\mathcal{M}_{g,r}(p_2);\mathbf{Z})=\mathbf{Z}/2\mathbf{Z}\text{\ or \,}\mathbf{Z}/4\mathbf{Z}.
\]
In particular if genus $g$ is even, 
\[
H_1(\mathcal{M}_g(p_2);\mathbf{Z})=\mathbf{Z}/2\mathbf{Z}.
\]
From the isomorphism $\hat{\mathcal{M}}_{(g,1)}(p_2)\cong\mathcal{M}_{g,1}(p_2)$ and the surjective homomorphism $b:\hat{\mathcal{M}}_{(g,1)}(p_2)\to\hat{\mathcal{M}}_{(g)}(p_2)$, we have 
\[
H_1(\hat{\mathcal{M}}_{(g,r)}(p_2);\mathbf{Z})=\mathbf{Z}/2\mathbf{Z}\text{\ or \,}\mathbf{Z}/4\mathbf{Z}
\]
for $r=0,1$.  
\begin{remark}\label{rem:generator}
For $r=0,1$, pick a simple closed curve $c\subset \Sigma_{g,r}$.
If the intersection number $c\cdot B_g$ is odd, then $[T_c^2]\in H_1(\Gamma_g(p_2);\mathbf{Z})$ is a generator. Hence $[T_c^2]\in H_1(\mathcal{M}_{g,r}(p_2);\mathbf{Z})$ is also a generator, and the lift of $T_c^2$ is a generator of $H_1(\hat{\mathcal{M}}_{(g,r)}(p_2);\mathbf{Z})$.

If $c$ is included in the subsurface $S$, we have $[T_c^2]=0\in H_1(\hat{\mathcal{M}}_{(g,r)}(p_2);\mathbf{Z})$, by Proposition \ref{generator}.
\end{remark}
\newpage
\section{A surjective homomorphism $\hat{\mathcal{M}}_{(g)}(p_2)\to\mathbf{Z}/4\mathbf{Z}$}\label{surj}
For a root of unity $\zeta$, we denote by $<\!\!\!\>\zeta\>\!\!\!>$ the cyclic group generated by $\zeta$.
In this section, we construct a surjective homomorphism 
\[
e:\hat{\mathcal{M}}_{(g)}(p_2)\to<\!\sqrt{-1}\!>
\]
using the Schottky theta constant associated with the cover $p_2:\Sigma_{2g-1}\to\Sigma_g$ when $g\ge2$, to complete Theorem \ref{main-theorem}. In the following, suppose genus $g\ge2$. 
\subsection{The Jacobi variety and the Prym variety}
Endow the surface $\Sigma_g$ with the structure of a Riemann surface $R$. Then the covering map $p_2:\Sigma_{2g-1}\to\Sigma_g$ induces the structure of a Riemann surface $\hat{R}$ in the surface $\Sigma_{2g-1}$. In this subsection, we review the Jacobi variety of the Riemann surface $R$ and the Prym variety of the double unbranched cover $p_2:\hat{R}\to R$. 
\begin{definition}
A $g$-characteristic is a row vector $m\in\mathbf{Z}^{2g}$. We denote $m=(m'|m'')$ where $m'=(m'_1,m'_2,\cdots,m'_g)$, $m''=(m''_1,m''_2,\cdots,m''_g)\in \mathbf{Z}^g$. We call the $g$-chatacteristic $m$ is even (resp. odd) if $\sum_{i=1}^g m'_im''_i$ is even (resp. odd).
\end{definition}
We denote the Siegel upper half space of degree $g$ by $\mathfrak{S}_g$. For a $g$-characteristic $m=(m'|m'')\in \mathbf{Z}^{2g}$ and $\tau\in \mathfrak{S}_g$, $z\in \mathbf{C}^g$, The theta function $\theta_m$ is defined by 
\[\theta_m(\tau, z):=\sum_{p\in\mathbf{Z}^{g}}\exp(\pi  i\{(p+m'/2)\tau\leftidx{^t}(p+m'/2)+(p+m'/2)\leftidx{^t}(z+m''/2)\}).\]
We denote $\theta_m(\tau, 0)$ simply by $\theta_m(\tau)$. Let $\Omega$ be the sheaf of holomorphic 1-forms on $R$. Choose a symplectic basis $\{A_i, B_i\}_{i=1}^g$ of $H_1(R;\mathbf{Z})$. It is known that under the homomorphism  
\[
\begin{array}{cccc}
H_1(R;\mathbf{Z})&\to &H^0(R;\Omega)^*&:=\Hom(H^0(R;\Omega),\mathbf{C}),\\
c&\mapsto&(\omega\mapsto\int_c\omega)
\end{array}
\]
$H_1(R;\mathbf{Z})$ maps onto a lattice in $H^0(R;\Omega)^*$. The Jacobi variety of $R$ is defined by
\[J(R)=\frac{H^0(R;\Omega)^*}{H_1(R;\mathbf{Z})}.\]
A basis $\{\omega_i\}_{i=1}^g$ of $H^0(R;\Omega)$ that satisfies
\[\int_{A_j}\omega_i=
\begin{cases}
1,&\text{ if }i=j,\\
0,&\text{ if }i\ne j,
\end{cases}
\]
is called the normalized basis with respect to the symplectic basis $\{A_i, B_i\}_{i=1}^g$. For the normalized basis $\{\omega_i\}_{i=1}^g$, the $g$-square matrix
\[\tau=(\tau_{ij}), \hspace{0.5cm} \tau_{ij}=\int_{B_j}\omega_i\]
is known to be the elements of the Siegel upper half space $\mathfrak{S}_g$, and is called the period matrix. For an even $g$-characteristic $m=(m'|m'')$ and the period matrix $\tau$, $\theta_{m}(\tau)$ is called the Riemann theta constant with $m$ associated with the compact Riemann surface $R$ and $\{A_i, B_i\}_{i=1}^g$.

Denote the generator of the deck transformation group of the cover $\hat{R}\to R$ by $t: \hat{R}\to \hat{R}$, the $(-1)$-eigenspace of $t_*: H_1(\hat{R};\mathbf{Z})\to H_1(\hat{R};\mathbf{Z})$ by
\[
H_1(\hat{R};\mathbf{Z})^{-}=\{c\in H_1(\hat{R};\mathbf{Z})\ |\ t_*(c)=-c \}, 
\]
and the $(-1)$-eigenspace of $t^*: H^0(\hat{R};\Omega)\to H^0(\hat{R};\Omega)$ by
\[
H^0(\hat{R};\Omega)^{-}=\{\omega\in H^0(\hat{R};\Omega)\ |\ t^*(\omega)=-\omega \}.
\]
Under the homomorphism  
\[
\begin{array}{cccc}
H_1(\hat{R};\mathbf{Z})&\to &H^0(\hat{R};\Omega)^*&:=\Hom(H^0(\hat{R};\Omega),\mathbf{C}),\\
c&\mapsto&(\omega\mapsto\int_c\omega)
\end{array}
\]
$H_1(\hat{R};\mathbf{Z})^-$ maps onto a lattice in $(H^0(R;\Omega)^-)^*$. 
\begin{definition}
The Prym variety $\Prym(\hat{R}, p_2)$ of the cover $p_2$ is defined by 
\[\Prym(\hat{R}, p_2)=\frac{(H^0(\hat{R};\Omega)^{-})^*}{H_1(\hat{R};\mathbf{Z})^{-}}\subset J(\hat{R}).\]
\end{definition}
For a symplectic basis $\{A_i, B_i\} _{i=1}^g$ of $H_1(R;\mathbf{Z})$, we choose a basis of $H_1(\hat{R};\mathbf{Z})$ as follows. For $i=1,2,\cdots,g-1$, denote the two lifts of $A_i$ by $\hat{A}_i$ and $\hat{A}_{i+g}$, and the two lifts of $B_i$ by  $\hat{B}_i$ and $\hat{B}_{i+g}$ such that 
\[
\hat{A}_i\cdot \hat{B}_i=1.
\]
The lifts of $2A_g$ and $B_g$ are uniquely determined, and denote them by $\hat{A}_g$ and $\hat{B}_g$, respectively. Then, $\{A_i-A_{g+i}, B_i-B_{g+i}\}_{i=1}^{g-1}$ form a basis of $H_1(\hat{R};\mathbf{Z})^{-}$. Moreover since the basis $\{\hat{A}_i-\hat{A}_{g+i}, \hat{B}_i-\hat{B}_{g+i}\}_{i=1}^{g-1}$ of $H_1(\hat{R};\mathbf{Z})^{-}$ satisfies
\begin{gather*}
(\hat{A}_i-\hat{A}_{g+i})\cdot(\hat{A}_j-\hat{A}_{g+j})=0,\ (\hat{B}_i-\hat{B}_{g+i})\cdot(\hat{B}_j-\hat{B}_{g+j})=0\\
(\hat{A}_i-\hat{A}_{g+i})\cdot(\hat{B}_j-\hat{B}_{g+j})=2\delta_{i\,j}.
\end{gather*}
Therefore, the action of $\hat{\varphi}\in\hat{\mathcal{M}}_{(g)}(p_2)$ on the basis $\{\hat{A}_i-\hat{A}_{g+i}, \hat{B}_i-\hat{B}_{g+i}\}_{i=1}^{g-1}$ induces the homomorphism
\[\tilde{\iota}:\hat{\mathcal{M}}_{(g)}(p_2)\to\Sp(2g-2;\mathbf{Z}).\]

For the above symplectic basis $\{\hat{A}_i, \hat{B}_i\}_{i=1}^{2g-1}$, choose the normalized basis $\{\hat{\omega}_i\}_{i=1}^{2g-1}$ of $H^0(\hat{R};\Omega)$, then $\{(\hat{\omega}_i-\hat{\omega}_{g+i})/2\}_{i=1}^{g-1}$ is a basis of $H^0(\hat{R};\Omega)^{-}$. It is known that the $(g-1)$-square matrix\[\tilde{\tau}=(\tilde{\tau}_{ij}), \hspace{0.5cm} \tilde{\tau}_{ij}=\int_{\hat{B}_j-\hat{B}_{g+j}}\frac{\hat{\omega}_i-\hat{\omega}_{g+i}}{2}\]
is the element of the Siegel upper half space $\mathfrak{S}_{g-1}$. We call $\tilde{\tau}$ the period matrix of the Prym variety.
\begin{definition}
For an even $(g-1)$-characteristic $\tilde{m}=(\tilde{m}'|\tilde{m}'')$ and the period matrix $\tilde{\tau}$ of $\Prym(\hat{R}, p_2)$, $\theta_{\tilde{m}}(\tilde{\tau})$ is called the Schottky theta constant with $\tilde{m}$ associated with the cover $p_2:\hat{R}\to R$ and $\{\hat{A}_i, \hat{B}_i\}_{i=1}^{2g-1}$.
\end{definition}
\subsection{Definition of the homomorphism $e:\hat{\mathcal{M}}_g(p_2)\to<\!\sqrt{-1}\!>$}
In this subsection, we give the definition of the homomorphism $e:\hat{\mathcal{M}}_{(g)}(p_2)\to<\!\sqrt{-1}\!>$.
Let $\tau$ be the period matrix of $R$, and $\tilde{\tau}$ the period matrix of the cover $p_2$. Consider the function
\[
\Phi_{m, n}^{\tilde{m}}(\tilde{\tau}, \tau)=\frac{\tilde\theta_{\tilde{m}}^2(\tilde{\tau})}{\theta_m(\tau)\theta_n(\tau)}
\]
for even g-characteristics $m,n$ and an even $(g-1)$-characteristic $\tilde{m}$. For a generic Riemann surface and a double unbranched covering space, $\Phi_{m,n}^{\tilde{m}}(\tilde{\tau}, \tau)$ is known to be a nonzero complex number (Fay\cite{fay1973tfr}). For a $g$-square matrix $M=(m_{ij})$, denote the row vector obtained by taking the diagonal entries of $M$ by $M_0:=(m_{11}, m_{22}, \cdots, m_{gg})\in\mathbf{Z}^g$. For 
$\sigma=
\begin{pmatrix}
\alpha&\beta\\
\gamma&\delta
\end{pmatrix}
\in\Sp(2g;\mathbf{Z})$ and a $g$-characteristic $m$,  
we define 
\[
\sigma\cdot m=m
\begin{pmatrix}
\leftidx{^t}{\alpha}{}&-\leftidx{^t}{\gamma}{}\\
-\leftidx{^t}{\beta}{}&\leftidx{^t}{\delta}{}
\end{pmatrix}
+((\leftidx{^t}{\beta}{}\alpha)_0|(\leftidx{^t}{\delta}{}\gamma)_0)\in \mathbf{Z}^{2g}.
\]
Note that this is not an action of $\Sp(2g;\mathbf{Z})$ on $\mathbf{Z}^{2g}$, and that this definition is different from that of Igusa\cite{igusa1964grt}.  
For $\hat{\varphi}\in\hat{\mathcal{M}}_{(g)}(p_2)$, denote $P_2(\hat{\varphi})$ by $\varphi\in\mathcal{M}_g(p_2)$. For an even $(g-1)$-characteristic $\tilde{m}$, choose the $g$-characteristics $m=(\tilde{m}',0|\tilde{m}'',1)$ and $n=(\tilde{m}',0 |\tilde{m}'',0)$.  Define the map $d_{\tilde{m},(\tilde{\tau}, \tau)}:\hat{\mathcal{M}}_{(g)}(p_2)\to\mathbf{C}$ by
\[
d_{\tilde{m}, (\tilde{\tau}, \tau)}(\hat{\varphi}):=
\frac{\Phi_{m,n}^{\tilde{m}}(\tilde{\tau}, \tau)}{\Phi_{\iota(\varphi)\cdot m,\iota(\varphi)\cdot n}^{\tilde{\iota}(\hat{\varphi})\cdot \tilde{m}}(\tilde{\tau}, \tau)}.
\]
In the next subsection, we will prove that $d_{\tilde{m}}=d_{\tilde{m}, (\tilde{\tau}, \tau)}$ is independent of the period matrices $\tilde{\tau}$ and $\tau$, and that the image of $d_{\tilde{m}}$ equals $<\!-1\!>$. For
$\sigma=
\begin{pmatrix}
\alpha&\beta\\
\gamma&\delta
\end{pmatrix}\in\Sp(2g;\mathbf{Z})$ and $\tau\in\mathfrak{S}_g$, we define the action of $\Sp(2g;\mathbf{Z})$ on $\mathfrak{S}_g$ by
\[
\sigma\cdot\tau:=(\delta\tau+\gamma)(\beta\tau+\alpha)^{-1}.
\] 
For $\sigma=\begin{pmatrix}\alpha&\beta\\\gamma&\delta\end{pmatrix}\in\Sp(2g;\mathbf{Z})$, it is known that the theta function has the transformation law (see Igusa\cite{igusa1972tf})
\[
\theta_{\sigma\cdot m}(\sigma\cdot\tau)=\gamma_m(\sigma)\det(\beta\tau+\alpha)^{-\frac{1}{2}}\theta_m(\tau),
\]
where $\gamma_m(\sigma)\in<\!\exp(\pi/4)\!>$ is called the theta multiplier. 

Now we can construct a homomorphism $e_{\tilde{m}}: \hat{\mathcal{M}}_{(g)}(p_2)\to <\!\sqrt{-1}\!>$ using $d_{\tilde{m}}$ and $\gamma_m$. For $\hat{\varphi}\in\hat{\mathcal{M}}_{(g)}(p_2)$ and an even $(g-1)$-characteristic $\tilde{m}$, define the map $e_{\tilde{m}}$ by
\[
e_{\tilde{m}}(\hat{\varphi}):=d_{\tilde{m}}(\hat{\varphi})\frac{\gamma_{\tilde{m}}^2(\tilde{\iota}(\hat{\varphi}))}{\gamma_m(\iota(\varphi))\gamma_n(\iota(\varphi))}.
\]
Note that $\gamma_m^2(\iota(\varphi))$ and $\gamma_m(\iota(\varphi))\gamma_n(\iota(\varphi))$ are uniquely determined. We will prove that $e=e_{\tilde{m}}$ is a homomorphism independent of the choice of $\tilde{m}$, and that the image of $e_{\tilde{m}}$ equals $<\!\sqrt{-1}\!>$ in the next subsection. 
\subsection{Proof of the main theorem}
In this subsection, we will prove that $e_{\tilde{m}}:\hat{\mathcal{M}}_{(g)}(p_2)\to <\!\sqrt{-1}\!>$ is a surjective homomorphism. We also prove that $d_{\tilde{m}}=d_{\tilde{m},(\tilde{\tau},\tau)}$ does not depends on the choice of $(\tilde{\tau},\tau)$, and that the image of $d_{\tilde{m}}$ equals the cyclic group $<\!-1\!>$.  For $\hat{\varphi}\in\mathcal{M}_{(g)}(p_2)$, we denote simply $\tilde{\iota}(\hat{\varphi})\in\Sp(2g-2;\mathbf{Z})$ and $\iota(\varphi)=\iota(P(\hat{\varphi}))\in\Gamma_g(p_2)$ by $\tilde{\sigma}$ and $\sigma$, respectively. 

To prove that $d_{\tilde{m}}$ only depends on $\hat{\varphi}\in\hat{\mathcal{M}}_{(g)}(p_2)$ and $\tilde{m}\in \mathbf{Z}^{g-1}$, we need the following theorem.
\begin{theorem}[Farkas, Rauch\cite{farkas1970prs}]\label{Farkas}\ \\
For an even $(g-1)$-characteristic $\tilde{m}$, define the $g$-characteristics $m=(\tilde{m}',0|\tilde{m}'',1)$ and $n=(\tilde{m}',0 |\tilde{m}'',0)$. Then  $\Phi_{m,n}^{\tilde{m}}(\tilde{\tau}, \tau)$ does not depend on the choice of $\tilde{m}$.
\end{theorem}
Define $\pi: \mathbf{Z}^{2g}\to\mathbf{Z}^{2g-2}$ by $\pi(m'|m'')=(m'_1,m'_2,\cdots, m'_{g-1}\ |\ m''_1, m''_2,\cdots, m''_{g-1})$. 
\begin{lemma}\label{lemma:mod2char}
For an even $g$-characteristic $\tilde{m}$ and $\hat{\varphi}\in\mathcal{M}_{(g)}(p_2)$,
\[\tilde{\sigma}\cdot \tilde{m}\equiv \pi(\sigma\cdot m)\ \mod 2,\]
where $m=(\tilde{m}',0|\tilde{m}'',1)$.
\end{lemma}
\begin{proof}
Denote the $1$-eigenspace of $H_1(\hat{R};\mathbf{Q})$ by $H_1(\hat{R};\mathbf{Q})^+$. Then 
\[
\{\hat{A}_i+\hat{A}_{g+i},\ \hat{B}_i+\hat{B}_{g+i}\}_{i=1}^{g-1}\cup \{\hat{A}_g,\ 2\hat{B}_g\}
\]
is a basis of $H_1(\hat{R};\mathbf{Q})^+$. The restriction of $p_2$
\[
H_1(\hat{R};\mathbf{Q})^+\to H_1(R;\mathbf{Q})
\]
maps the basis $\{\hat{A}_i+\hat{A}_{g+i},\ \hat{B}_i+\hat{B}_{g+i}\}_{i=1}^{g-1}\cup \{\hat{A}_g,\ 2\hat{B}_g\}\in H_1(\hat{R};\mathbf{Q})^+$ to the basis $\{2A_i,\ 2B_i\}_{i=1}^{g}\in H_1(R;\mathbf{Q})$. Since for $i=1,\cdots,g-1$ we have
\begin{gather*}
\varphi_*(2A_i)=\varphi_*(p_2)_*(\hat{A}_i+\hat{A}_{g+i})= (p_2)_*\hat{\varphi}_*(\hat{A}_i+\hat{A}_{g+i}),\ \varphi(2A_g)= (p_2)_*\hat{\varphi}_*(\hat{A}_g),\\
\varphi_*(2B_i)= (p_2)_*\hat{\varphi}_*(\hat{B}_i+\hat{B}_{g+i}),\text{ and }\varphi_*(2B_g)= (p_2)_*\hat{\varphi}_*(2\hat{B}_g).
\end{gather*}
Hence, the induced homomorphism $\hat{\mathcal{M}}_{(g)}(p_2)\to\Sp(2g;\mathbf{Z})$ by the action of $\hat{\varphi}\in\hat{\mathcal{M}}_{(g)}(p_2)$ on the basis $\{\hat{A}_i+\hat{A}_{g+i},\ \hat{B}_i+\hat{B}_{g+i}\}_{i=1}^{g-1}\cup \{\hat{A}_g,\ 2\hat{B}_g\}$ is equal to $\iota P_2:\hat{\mathcal{M}}_{(g)}(p_2)\to\Gamma_g(p_2)$. Denote $\tilde{\sigma}\in\Sp(2g-2;\mathbf{Z})$ by
\[
\tilde{\sigma}=
\begin{pmatrix}
\alpha'&\beta'\\
\gamma'&\delta'
\end{pmatrix},
\]
where $\alpha', \beta', \gamma', \delta'\in M(g-1;\mathbf{Z})$. Since we have
\[
\hat{\varphi}_*(\hat{A}_i+\hat{A}_{g+i})\equiv\hat{\varphi}_*(\hat{A}_i-\hat{A}_{g+i}),\text{ and }\hat{\varphi}_*(\hat{B}_i+\hat{B}_{g+i})\equiv\hat{\varphi}_*(\hat{B}_i-\hat{B}_{g+i}) \ \mod 2,
\]
and $\sigma=\iota P_2(\hat{\varphi})$ preserves the homology class $B_g \mod2$, $\sigma$ is written in the form
\[
\sigma=
\begin{pmatrix}
\alpha'&\leftidx{^t}{v_1}{}&\beta'&0\\
0&1&0&0\\
\gamma'&\leftidx{^t}{v_2}{}&\delta'&0\\
v_3&k&v_4&1
\end{pmatrix}
\mod2,
\]
where $v_1,v_2,v_3,v_4\in \mathbf{Z}^{g-1}, k\in\mathbf{Z}$. Then it is easy to see that $\pi(\sigma\cdot m)\equiv\tilde{\sigma}\cdot \tilde{m}\ \mod2$.
\end{proof}
\begin{lemma}
For $\hat{\varphi}\in\hat{\mathcal{M}}_{(g)}(p_2)$, the value $d_{\tilde{m}}(\hat{\varphi})=d_{\tilde{m},(\tilde{\tau},\tau)}(\hat{\varphi})$ does not depend on the choice of $(\tilde{\tau},\tau)$, and the image of $d_{\tilde{m}}$ equals the cyclic group $<\!-1\!>$. In particular, it does not depend on a complex structure of the cover $p_2:\hat{R}\to R$.
\end{lemma}
\begin{proof}
Note that, for any g-characteristic $u=(u'|u''), v=(v'|v'')\in\mathbf{Z}$ we have
\[\theta_{u+2v}=(-1)^{u'v''}\theta_u,\]
by the definition of the theta function. Consider the g-characteristic $v_0=(0,\cdots,0,1|0,\cdots,0,0)\in\mathbf{Z}^{2g}$. Since $\sigma$ preserves the homology class $B_g\mod2$, we have
\begin{gather*}
\sigma\cdot(m-n)=(m-n)
\begin{pmatrix}
\leftidx{^t}{\alpha}{}&-\leftidx{^t}{\gamma}{}\\
-\leftidx{^t}{\beta}{}&\leftidx{^t}{\delta}{}
\end{pmatrix}
\equiv v_0
\begin{pmatrix}
\leftidx{^t}{\alpha}{}&-\leftidx{^t}{\gamma}{}\\
-\leftidx{^t}{\beta}{}&\leftidx{^t}{\delta}{}
\end{pmatrix}
\equiv v_0\ {, and}\\
(\sigma\cdot m)'_g\equiv(\sigma\cdot n)'_g\equiv(\beta\leftidx{^t}{\alpha}{})_{gg}\equiv0\ \mod 2.
\end{gather*}
By Lemma \ref{lemma:mod2char}, there exists $v_1, v_2\in\mathbf{Z}^{2g}$ such that 
\[
\sigma\cdot m+2v_1=((\tilde{\sigma}\cdot \tilde{m})',0|(\tilde{\sigma}\cdot \tilde{m})'',k_1),\text{ and } \sigma\cdot n+2v_2=((\tilde{\sigma}\cdot \tilde{m})',0|(\tilde{\sigma}\cdot \tilde{m})'',k_2),
\]
where 
\[
k_1=0\text{ or }1,\ k_2=0\text{ or }1,\text{ and } k_1+k_2=1.
\]
Then there exists $p(\tilde{m},\hat{\varphi})\in <\!-1\!>$ such that
\[
\Phi_{\sigma\cdot m+2v_1, \sigma\cdot n+2v_2}^{\tilde{\sigma}\cdot\tilde{m}}(\tilde{\tau}, \tau)=p(\tilde{m}, \hat{\varphi})\Phi_{\sigma\cdot m, \sigma\cdot n}^{\tilde{\sigma}\cdot\tilde{m}}(\tilde{\tau}, \tau).
\]
Note that $p(\tilde{m},\hat{\varphi})$ does not depend on the choice of $(\tilde{\tau},\tau)$. By Theorem \ref{Farkas}, we have
\[
\Phi_{\sigma\cdot m+2v_1, \sigma\cdot n+2v_2}^{\tilde{\sigma}\cdot\tilde{m}}(\tilde{\tau}, \tau)=\Phi_{m,n}^{\tilde{m}}(\tilde{\tau}, \tau).
\]
Hence we have
\[
p(\tilde{m},\hat{\varphi})=d_{\tilde{m}}(\hat{\varphi}).
\]
This proves the lemma.
\end{proof}
Consider the action of $\varphi\in\mathcal{M}_g(p_2)$ on the symplectic basis $\{A_i, B_i\}_{i=1}^g$. The basis $\{\varphi_*A_i, \varphi_*B_i\} _{i=1}^g$ is also a symplectic basis of $H_1(R;\mathbf{Z})$. The corresponding period matrix is
\[
\tau'=(\tau'_{ij}), \hspace{0.5cm} \tau'_{ij}=\int_{\varphi_*B_j}\omega'_i,
\]
where $\{\omega'_i\}_{i=1}^{g}$ is the normalized basis. This is equal to $\leftidx{^t}{\iota(\varphi)}{}\cdot\tau$. Next, Consider the action of $\hat{\varphi}\in\hat{\mathcal{M}}_{(g)}(p_2)$ on the basis $\{\hat{A}_i, \hat{B}_i\} _{i=1}^{2g-1}$ of $H_1(\hat{R};\mathbf{Z})$. Note that the basis $\{\hat{\varphi}_*\hat{A}_i,\hat{\varphi}_*\hat{B}_i\}_{i=1}^{2g-1}$ is again the lift of $\{\varphi_*A_i, \varphi_*B_i\} _{i=1}^g$. The period matrix of $\Prym(\hat{R}, p_2)$ with respect to the basis $\{\hat{\varphi}_*(\hat{A}_i-\hat{A}_{g+i}), \hat{\varphi}_*(\hat{B}_i-\hat{B}_{g+i})\} _{i=1}^{2g-1}$ of $H_1(\hat{R};\mathbf{Z})$ is
\[
\tilde{\tau}':=(\tilde{\tau}'_{ij}), \hspace{0.5cm} \tilde{\tau}'_{ij}=\int_{\hat{\varphi}_*(\hat{B}_j-\hat{B}_{g+j})}\frac{\hat{\omega'}_i-\hat{\omega'}_{g+i}}{2},
\]
where $\{\hat{\omega'}_i\}_{i=1}^{2g-1}$ is the normalized basis. This is equal to $\leftidx{^t}{\tilde{\iota}}{}(\hat{\varphi})\cdot\tilde{\tau}$. Hence, $\leftidx{^t}{\iota(\varphi)}{}\cdot\tau$ is also the perod matrix of $R$, and $\leftidx{^t}{\tilde{\iota}}{}(\hat{\varphi})\cdot\tilde{\tau}$ is also the period matrix of the cover $p_2$. This shows that the pair $(\tilde{\sigma}\cdot\tilde{\tau}, \sigma\cdot\tau)$ satisfies the condition of Theorem \ref{Farkas} for any $\hat{\varphi}\in\hat{\mathcal{M}}_{(g)}(p_2)$. 
\begin{theorem}
The map $e_{\tilde{m}}$ is a homomorphism, and the image of $e_{\tilde{m}}(\hat{\varphi})$ equals $<\!\sqrt{-1}\!>$. Moreover $e({\hat{\varphi}}):=e_{\tilde{m}}(\hat{\varphi})$ does not depend on the choice of $\tilde{m}$.
\end{theorem}
\begin{proof}
For $\hat{\varphi}\in\hat{\mathcal{M}}_{(g)}(p_2)$, denote $\sigma_1:=\sigma=\iota P_2(\hat{\varphi}')$, and $\tilde{\sigma}_1:=\tilde{\sigma}=\tilde{\iota}(\hat{\varphi})$. Similarly, denote $\sigma_2:=\iota P_2(\hat{\varphi}')$, $\tilde{\sigma}_2:=\tilde{\iota}(\hat{\varphi}')$, and $\sigma_3:=\iota P_2(\hat{\varphi}\hat{\varphi}')$, $\tilde{\sigma}_3:=\tilde{\iota}(\hat{\varphi}\hat{\varphi}')$. Write $\sigma_i$ as
\[
\sigma_i=
\begin{pmatrix}
\alpha_i&\beta_i\\
\gamma_i&\delta_i
\end{pmatrix}
\hspace{0.5cm}
\text{for } i=1,2,3.
\]
We also denote simply $\tilde{\tau}':=\tilde{\sigma}_2\cdot\tilde{\tau}$, and
$\tau':=\sigma_2\cdot\tau$.
Since the pairs $(\tilde{\sigma}_1\tilde{\sigma}_2\cdot\tilde{\tau}, \sigma_1\sigma_2\cdot\tau)$, and $(\tilde{\sigma}_2\cdot\tilde{\tau}, \sigma_2\cdot\tau)$ satisfies the condition of Theorem \ref{Farkas}, we have 
\begin{gather*}
\frac{1}{d_{\tilde{m}}(\hat{\varphi}\hat{\varphi'})}
=\frac{\Phi_{(\sigma_1\sigma_2)\cdot m, (\sigma_1\sigma_2)\cdot n}^{(\tilde{\sigma}_1\tilde{\sigma}_2)\cdot\tilde{m}}(\tilde{\sigma}_1\tilde{\sigma}_2\cdot\tilde{\tau}, \sigma_1\sigma_2\cdot\tau)}{\Phi_{m, n}^{\tilde{m}}(\tilde{\sigma}_1\tilde{\sigma}_2\cdot\tilde{\tau}, \sigma_1\sigma_2\cdot\tau)}
=\frac{\gamma_{\tilde{m}}^2(\tilde{\sigma}_1\tilde{\sigma}_2)}{\gamma_m(\sigma_1\sigma_2)\gamma_n(\sigma_1\sigma_2)}
\frac{\det(\tilde{\beta}_3\tilde{\tau}+\tilde{\alpha}_3)^{-1}}{\det(\beta_3\tau+\alpha_3)^{-1}}
\frac{\Phi_{m, n}^{\tilde{m}}(\tilde{\tau}, \tau)}{\Phi_{m, n}^{\tilde{m}}(\tilde{\sigma}_1\tilde{\sigma}_2\cdot\tilde{\tau}, \sigma_1\sigma_2\cdot\tau)},\\
\frac{1}{d_{\tilde{m}}(\hat{\varphi})}
=\frac{\Phi_{\sigma_1\cdot m, \sigma_1\cdot n}^{\tilde{\sigma}_1\cdot\tilde{m}}(\tilde{\sigma}_1\tilde{\sigma}_2\cdot\tilde{\tau}, \sigma_1\sigma_2\cdot\tau)}{\Phi_{m, n}^{\tilde{m}}(\tilde{\sigma}_1\tilde{\sigma}_2\cdot\tilde{\tau}, \sigma_1\sigma_2\cdot\tau)}
=\frac{\gamma_{\tilde{m}}^2(\tilde{\sigma}_1)}{\gamma_m(\sigma_1)\gamma_n(\sigma_1)}
\frac{\det(\tilde{\beta}_1\tilde{\tau}'+\tilde{\alpha}_1)^{-1}}{\det(\beta_1\tau'+\alpha_1)^{-1}}
\frac{\Phi_{m, n}^{\tilde{m}}(\tilde{\sigma}_2\cdot\tilde{\tau}, \sigma_2\cdot\tau)}{\Phi_{m, n}^{\tilde{m}}(\tilde{\sigma}_1\tilde{\sigma}_2\cdot\tilde{\tau}, \sigma_1\sigma_2\cdot\tau)},\\
\frac{1}{d_{\tilde{m}}(\hat{\varphi'})}
=\frac{\Phi_{\sigma_2\cdot m, \sigma_2\cdot n}^{\tilde{\sigma}_2\cdot\tilde{m}}(\tilde{\sigma}_2\cdot\tilde{\tau}, \sigma_2\cdot\tau)}{\Phi_{m, n}^{\tilde{m}}(\tilde{\sigma}_2\cdot\tilde{\tau}, \sigma_2\cdot\tau)}
=\frac{\gamma_{\tilde{m}}^2(\tilde{\sigma}_2)}{\gamma_m(\sigma_2)\gamma_n(\sigma_2)}
\frac{\det(\tilde{\beta}_2\tilde{\tau}+\tilde{\alpha}_2)^{-1}}{\det(\beta_2\tau+\alpha_2)^{-1}}
\frac{\Phi_{m, n}^{\tilde{m}}(\tilde{\tau}, \tau)}{\Phi_{m, n}^{\tilde{m}}(\tilde{\sigma}_2\cdot\tilde{\tau}, \sigma_2\cdot\tau)},
\end{gather*}
by the definition of $d_{\tilde{m}}(\hat{\varphi})$.
It is easy to see that
\begin{gather*}
\det(\tilde{\beta}_2\tilde{\tau}+\tilde{\alpha}_2)\det(\tilde{\beta}_1\tilde{\tau}'+\tilde{\alpha}_1)=\det(\tilde{\beta}_3\tilde{\tau}+\tilde{\alpha}_3)\text{, and}\\
\det(\beta_2\tau+\alpha_2)\det(\beta_1\tau'+\alpha_1)=\det(\beta_3\tau+\alpha_3).
\end{gather*}
This shows that $e_{\tilde{m}}$ is a homomorphism.

Next, we determine the image of $e_{\tilde{m}}$. There are two lifts in $\hat{\mathcal{M}}_{(g)}(p_2)$ of $a^2\in\mathcal{M}_g(p_2)$. We denote the lift which fix the homology class $\hat{A_1}$ by $\hat{a}\in\hat{\mathcal{M}}_{(g)}(p_2)$. As we stated in Remark \ref{rem:generator}, $H_1(\hat{\mathcal{M}}_{(g)}(p_2);\mathbf{Z})$ is generated by $\hat{a}$. For $\hat{\varphi}=\hat{a}$, we have $\tilde{\sigma}=\tilde{\iota}(\hat{a})=I_{2g-2}\in \Sp(2g-2;\mathbf{Z})$, $\sigma=\iota P_2(\hat{a})=\gamma_{gg}\in\Gamma_g(p_2)$. From Theorem 3 in Igusa\cite{igusa1964grt}, for any $\tilde{m}\in\mathbf{Z}^{2(g-1)}$, we have
\[
\gamma_m(\sigma)\gamma_n(\sigma)=-\sqrt{-1}, \text{ and \ } \gamma_{\tilde{m}}^2(\tilde{\sigma})=1,
\]
so that 
\[
\frac{\gamma_{\tilde{m}}^2(\tilde{\sigma})}{\gamma_m(\sigma)\gamma_n(\sigma)}=\sqrt{-1}.
\]
It is easy to see that $d_{\tilde{m}}(\hat{a})=1$. Hence $e_{\tilde{m}}(\hat{a})$ is a generator of the cyclic group $<\!\sqrt{-1}\!>$ and is independent of the choice of $\tilde{m}$. 
\end{proof}
For $r=0,1$, we proved $H_1(\hat{\mathcal{M}}_{(g,r)}(p_2);\mathbf{Z})\cong\mathbf{Z}/2\mathbf{Z}\text{ or }\mathbf{Z}/4\mathbf{Z}$ in Section \ref{genGamma}. From the above Theorem, we have 
\[H_1(\hat{\mathcal{M}}_{(g,r)}(p_2);\mathbf{Z})\cong\mathbf{Z}/4\mathbf{Z}.\]
Since $\mathcal{M}_{g,1}(p_2)$ is isomorphic to $\hat{\mathcal{M}}_{(g,1)}(p_2)$, we have $H_1(\mathcal{M}_{g,1}(p_2);\mathbf{Z})\cong \mathbf{Z}/4\mathbf{Z}$. Consider $H_1(\mathcal{M}_g(p_2);\mathbf{Z})$ when genus $g$ is odd. For the deck transformation $t$, we obtain
\[
e(t)=(-1)^{g-1},
\]
from Theorem 3 in Igusa\cite{igusa1964grt}. By the Lyndon-Hochschild-Serre spectral sequence, we have
\[
\mathbf{Z}/2\mathbf{Z}\to H_1(\hat{\mathcal{M}}_{(g)}(p_2);\mathbf{Z})\to H_1(\mathcal{M}_g(p_2);\mathbf{Z})\to 0
\]
This shows that $H_1(\mathcal{M}_g(p_2);\mathbf{Z})\cong\mathbf{Z}/4\mathbf{Z}$ when $g$ is odd. This completes the proof of Theorem \ref{main-theorem}. 
From the Theorem \ref{main-theorem}, we obtain many homomorphisms $\mathcal{M}_{g,1}[d]\to\mathbf{Z}/4\mathbf{Z}$ for an even integer $d$. 
\begin{proposition}\label{prop:leveld}
For a positive even integer $d$, there exists an injection 
\[
(\mathbf{Z}/4\mathbf{Z})^{2g} \hookrightarrow \Hom(\mathcal{M}_{g,1}[d];\mathbf{Z}/4\mathbf{Z}).
\]
When $d=2$ and $g$ is $odd$, we have
\[
(\mathbf{Z}/4\mathbf{Z})^{2g} \hookrightarrow \Hom(\mathcal{M}_{g}[d];\mathbf{Z}/4\mathbf{Z}).
\]
\end{proposition}
\begin{proof}
To prove the proposition, we will construct a homomorphism from $\mathcal{M}_{g,1}[d]$ into $\mathcal{M}_{dg/2-1,1}(p'_X)$ for a certain double cover $p'_X$.

Let $X$ be one of the homology classes $A_1,\cdots,A_g,B_1,\cdots,B_g\in H_1(\Sigma_g;\mathbf{Z})$. Consider the $d$ cover $q_X:\Sigma_{dg-1}\to\Sigma_{g}$ such that the monodromy homomorphisms $\pi_1(\Sigma_{g})\to\mathbf{Z}/d\mathbf{Z}$ is equal to the Poincar\'e dual of $X\in H^1(\Sigma_{g};\mathbf{Z}/d\mathbf{Z})$. Denote a generator of the deck transformation group by $t_X$. Consider 
\[
\Sigma_{g,1}=\Sigma_g-D^2\subset \Sigma_g\text{, and }\Sigma_{dg-1,d}=\Sigma_{dg-1}-{q}_X^{-1}(D^2).
\]
We denote the restriction of the cover $q_X|_{\Sigma_{g,1}}:\Sigma_{dg-1,d}\to\Sigma_{g,1}$ by $p_X$.  Choose two connected components $D_1$ and $D_2$ of ${q}_X^{-1}(D^2)$ such that $t_X^{d/2}D_1=D_2$. Consider $\Sigma_{dg-1,2}=\Sigma_{dg-1}-\amalg_{i=1}^{2}D_i$. Then we have the double cover
\[p'_X:\Sigma_{dg-1,2}\to\Sigma_{dg-1,2}/<\!t_X^{d/2}\!>=\Sigma_{dg/2-1,1}.\]

We have the projection $P_X:\hat{\mathcal{M}}_{(g,1)}(p_X)\to\mathcal{M}_{g,1}(p_X)$ and $P'_X:\hat{\mathcal{M}}_{(dg/2-1,1)}(p'_X)\to\mathcal{M}_{dg/2-1,1}(p'_X)$. Since the centralizer of $<\!t_X\!>$ is included in the centralizer of $<\!t_X^{d/2}\!>$, we have the homomorphism 
\[
\begin{array}{cccc}
Q_X:&\hat{\mathcal{M}}_{(g,1)}(p_X)&\to&\hat{\mathcal{M}}_{(dg/2-1,1)}(p'_X).\\
&[\hat{f}]&\mapsto&[\hat{f}\cup id_{\cup_{i=1}^{d-2} D^2}]
\end{array}
\]
Note that we have the inclusion map $i_X:\mathcal{M}_{g,1}[d]\to\mathcal{M}_{g,1}(p_X)$. Hence we have the homomorphism
\[
P'_XQ_XP_X^{-1}i_X:\mathcal{M}_{g,1}[d]\to\mathcal{M}_{dg/2-1,1}(p'_X).
\]
Consider the induced homomorphism $(P'_XQ_XP_X^{-1}i_X)_*:H_1(\mathcal{M}_{g,1}[d];\mathbf{Z})\to H_1((\mathcal{M}_{dg/2-1,1}(p'_X);\mathbf{Z})$. For the simple closed curves $Y=A_1,\cdots,A_g,B_1,\cdots,B_g$, denote the Dehn twists along $Y$ by $T_Y$. Then we have
\[
(P'_XQ_XP_X^{-1}i_X)_*(T_Y^d)=
\begin{cases}
1,&\text{ if } Y=X,\\
0,&\text{ otherwise},
\end{cases}
\]
by Remark \ref{rem:generator}. Hence the induced map 
\[(\mathbf{Z}/4\mathbf{Z})^{2g} \to \Hom(\mathcal{M}_{g,1}[d];\mathbf{Z}/4\mathbf{Z})\]
is injective.

Next, consider the case of $d=2$ and $g$ is odd. Then $H_1(\mathcal{M}_{g}(p_X);\mathbf{Z})$ is isomorphic to $\mathbf{Z}/4\mathbf{Z}$. The inclusion $\mathcal{M}_{g}[2]\to\mathcal{M}_g(p_X)$ induces a homomorphism $H_1(\mathcal{M}_{g}[2];\mathbf{Z})\to H_1(\mathcal{M}_{g}(p_X);\mathbf{Z})\cong\mathbf{Z}/4\mathbf{Z}$. Similarly, we have the injective homomorphism $(\mathbf{Z}/4\mathbf{Z})^{2g} \to \Hom(\mathcal{M}_{g}[2];\mathbf{Z}/4\mathbf{Z})$. This completes the proof.
\end{proof}
\bibliographystyle{amsplain}
\bibliography{first.bib}
\end{document}